# Integral-Input-Output to State Stability

Brian Ingalls*
Dept. of Mathematics, Rutgers University, NJ

Department of Mathematics - Hill Center
Rutgers, The State University Of New Jersey
110 Frelinghuysen Rd.
Piscataway, NJ 08854-8019
FAX: 732-445-5530
ingalls@math.rutgers.edu

October 28, 2018

**Abstract**

A notion of detectability for nonlinear systems is discussed. Within the framework of "input to state stability" (ISS), a dual notion of "output to state stability" (OSS), and a more complete detectability notion, "input-output to state stability" (IOSS) have appeared in the literature. This note addresses a variant of the IOSS property, using an integral norm to measure signals, as opposed to the standard supremum norm that appears in ISS theory.

## 1 Introduction

We consider stability features for the system with output:

$$\dot{x}(t) = f(x(t), u(t)), \quad y(t) = h(x(t)), \tag{1}$$

where $x \in \mathbb{R}^n$. The function $f : \mathbb{R}^n \times \mathbb{R}^m \to \mathbb{R}^n$ is assumed jointly continuous in $x$ and $u$, and locally Lipschitz in $x$, uniformly in $u$. The output map $h : \mathbb{R}^n \to \mathbb{R}^p$ is assumed locally Lipschitz, and we suppose $f(0,0) = 0$ and $h(0) = 0$. Inputs $u(\cdot)$ take values in some set $\mathbb{U} \subseteq \mathbb{R}^m$ (where $\mathbb{U} = \mathbb{R}^m$ unless otherwise stated).

The notion of *input to state stability* (ISS), introduced in [22], provides a theoretical framework in which to formulate questions of robustness with respect to inputs (seen as disturbances) acting on a system. An ISS system is, roughly, one which has a "finite nonlinear gain" with respect to inputs and whose transient behavior can be bounded in terms of the size of the initial state; the precise definition is in terms of $\mathcal{K}$ function gains. The theory of ISS systems now forms an integral part of several texts ([4, 6, 8, 12, 13, 21]) as well as expository and research articles (see e.g. [7, 9, 14, 17] as well as the recent [27]).

*Supported in part by US Air Force Grant F49620-98-1-0242

Within the framework of ISS, a natural notion of detectability, and a dual to the ISS property, is the notion of *output to state stability* (OSS) addressed in [23, 24]. These references include a characterization of OSS in terms of a Lyapunov (or "storage") function, as well as a discussion of the roles of OSS and the more general property of input-output to state stability (IOSS) in nonlinear observer theory. The IOSS property was further addressed in [11].

In each of the notions mentioned thus far, signals (i.e. inputs and outputs) are measured by a supremum (or $L^\infty$) norm. In many cases, it may be more natural to use an integral (or $L^1$-type) norm, which corresponds to a measure of the "total energy" of the signal. A variant of ISS using this norm, called *integral-ISS* (iISS) was introduced in [26] and further studied in [1].

This paper addresses a combination of the ideas described above: namely a notion of detectability making use of integral norms. This property is formulated as a natural combination of the IOSS and iISS properties. It was introduced as integral-input-output to state stability (iIOSS) in [11]. This notion has been called "integral-detectability" by Morse and Hespanha [19] and is closely related to the notion of a "convergent observer" used by Krener in [10]. In addition, all systems which are passive in the sense of [13] automatically satisfy the iIOSS property (cf. remark 12 in [24]).

The main result in this paper is a characterization of the iIOSS property in terms of the existence of an appropriate Lyapunov function. In general, the result provides for the existence of a continuous Lyapunov function, though we indicate an important case where the construction can be extended to show the existence of a smooth function. Such Lyapunov characterizations for detectability notions are especially insightful, since in some cases the notion of detectability has been *defined* in terms of the existence of an appropriate Lyapunov (or "storage") function (e.g. [16, 18]).

While we refer to IOSS and iIOSS as notions of detectability, they should be called more precisely notions of *zero-detectability*, as they characterize the property that the information from the output is sufficient to deduce stability of the state to the origin. For linear systems, such a property is equivalent to "full-state detectability" – the property which allows construction of an observer which tracks arbitrary trajectories. For nonlinear systems, a zero-detectability condition cannot guarantee the existence of a "complete" observer. Given a nonlinear system which satisfies a zero-detectability property, the most one may expect is to be able to construct a *norm-observer* which is able to provide a bound on how far the state is from the origin. The existence of norm observers for IOSS systems was addressed in [11]. We shall see that a similar construction for iIOSS systems follows immediately from the definitions.

## 1.1 Basic Definitions and Notation

The Euclidean norm in a space $\mathbb{R}^k$ is denoted simply by $|\cdot|$. For each interval $\mathcal{I} \subseteq \mathbb{R}$ and any measurable function $u : \mathcal{I} \to \mathbb{R}^k$, we will use $\|u\|_{\mathcal{I}}$ to denote the (essential) supremum norm of $u(\cdot)$ over $\mathcal{I}$. That is, $\|u\|_{\mathcal{I}} = \operatorname{ess\,sup} \{|u(t)| \, : \, t \in \mathcal{I}\}$. An *input* (or *control*) will be a measurable, locally essentially bounded function $u : \mathcal{I} \to \mathbb{R}^m$, where $\mathcal{I}$ is a subinterval of $\mathbb{R}$ which contains the origin, such that $u(t) \in \mathbb{U} \subseteq \mathbb{R}^m$ for almost all $t \in \mathcal{I}$. Unless otherwise specified, we assume $\mathcal{I} = \mathbb{R}_{\geq 0}$.

For each initial state $\xi$ and input $u$ we let $x(t, \xi, u)$ denote the unique maximal solution of (1), and we write the output signal as $y(t, \xi, u) := h(x(t, \xi, u))$. A system is *forward complete* if each $\xi \in \mathbb{R}^n$ and each input $u$ defined on $\mathbb{R}_{\geq 0}$ produce a solution $x(t, \xi, u)$ which is defined for all $t \geq 0$.

A function $\gamma : \mathbb{R}_{\geq 0} \to \mathbb{R}_{\geq 0}$ is *of class* $\mathcal{K}$ (or a "$\mathcal{K}$ function") if it is continuous, positive definite, and strictly increasing; and is *of class* $\mathcal{K}_\infty$ if in addition it is unbounded. A function $\rho : \mathbb{R}_{\geq 0} \to \mathbb{R}_{\geq 0}$ is of class $\mathcal{L}$ if it is continuous, decreasing, and tends to zero as its argument tends to $+\infty$. A function $\beta : \mathbb{R}_{\geq 0} \times \mathbb{R}_{\geq 0} \to \mathbb{R}_{\geq 0}$ is *of class* $\mathcal{KL}$ if for each fixed $t \geq 0$, $\beta(\cdot, t)$ is of class $\mathcal{K}$ and for each fixed $s \geq 0$, $\beta(s, \cdot)$ is of class $\mathcal{L}$.

To formulate the statement that a nonsmooth function decreases in an appropriate manner, we will make use of the notion of the *viscosity subgradient* (cf. [3]).

**Definition 1.1** A vector $\zeta \in \mathbb{R}^n$ is a *viscosity subgradient* of the function $V : \mathbb{R}^n \to \mathbb{R}$ at $\xi \in \mathbb{R}^n$ if there exists a function $g : \mathbb{R}^n \to \mathbb{R}$ satisfying $\lim_{h \to 0} \frac{g(h)}{|h|} = 0$ and a neighbourhood $\mathcal{O} \subset \mathbb{R}^n$ of the origin so that

$$V(\xi + h) - V(\xi) - \zeta \cdot h \geq g(h)$$

for all $h \in \mathcal{O}$. □

The (possibly empty) set of viscosity subgradients of $V$ at $\xi$ is called the *viscosity subdifferential* and is denoted $\partial_D V(\xi)$. We remark that if $V$ is differentiable at $\xi$, then $\partial_D V(\xi) = \{\nabla V(\xi)\}$.

## 2 The integral-Input-Output to State Stability Property

The main property of interest in this paper is the following.

**Definition 2.1** We say that a forward complete system (1) satisfies the *integral-input-output to state stability* property (iIOSS) if there exist $\alpha \in \mathcal{K}_\infty$, $\beta \in \mathcal{KL}$, and $\gamma_1$, $\gamma_2 \in \mathcal{K}$ so that for every initial point $\xi \in \mathbb{R}^n$, and every input $u$,

$$\alpha(|x(t, \xi, u)|) \leq \beta(|\xi|, t) + \int_0^t \gamma_1(|y(s, \xi, u)|) + \gamma_2(|u(s)|)\, ds \tag{2}$$

for all $t \geq 0$. □

**Remark 2.2** We note that, by causality, the iIOSS bound (2) can be expressed *equivalently* as

$$\alpha(|x(t, \xi, u)|) \leq \beta(|\xi|, t) + \int_0^t \gamma_1(|y(s, \xi, u)|)\, ds + \int_0^\infty \gamma_2(|u(s)|)\, ds \tag{3}$$

for all $t \geq 0$. We will make use of this alternate description. □

**Remark 2.3** Recall that a forward complete system (1) satisfies the *input-output to state stability* property (IOSS) if there exist $\beta \in \mathcal{KL}$, and $\gamma_1$, $\gamma_2 \in \mathcal{K}$ so that for every initial point $\xi \in \mathbb{R}^n$, and every input $u$,

$$|x(t, \xi, u)| \leq \beta(|\xi|, t) + \gamma_1(\|y(\cdot, \xi, u)\|_{[0,t]}) + \gamma_2(\|u\|_{[0,t]})$$

for all $t \geq 0$.

It is natural to compare the notion of iIOSS to this analogous property. We will show as a consequence of our main result that an IOSS system is in particular an iIOSS system. It has been shown (in [11] and [1] respectively), that the iOSS property is strictly weaker than OSS, and that the iISS property is strictly weaker than ISS. Either of these results show that iIOSS is a *strictly* weaker property than IOSS. □

**Remark 2.4** It is an easy exercise to show that for linear systems the iIOSS property is equivalent to detectability. However, as mentioned above, for general systems as in (1), iIOSS is a notion of zero-detectability. Given that a system satisfies the iIOSS property, one cannot hope to build a complete observer for the system, but rather only a *norm observer* which measures how far the state is from the origin. The construction of norm observers for IOSS systems was addressed in [11], where it was shown that a system satisfies the IOSS property if and only if it admits an appropriate norm observer.

For iIOSS systems, the situation is more transparent. It is immediate that if a system satisfies the iIOSS bound (2), then for the function $p : \mathbb{R}_{\geq 0} \to \mathbb{R}_{\geq 0}$ defined by

$$\dot{p}(t) = \gamma_1(|y(t,\xi,u)|) + \gamma_2(|u(t)|), \qquad p(0) = 0,$$

the system will satisfy

$$\alpha(|x(t,\xi,u)|) \leq \beta(|\xi|\,,t) + p(t) \qquad \forall t \geq 0.$$

Thus the function $p(\cdot)$ provides an asymptotic upper bound for the size of the state, i.e. it is a norm observer for the system.

To extend these ideas to the case where construction of a full-state observer may be possible, one must consider a notion of "complete" detectability for nonlinear systems. Such a notion was introduced in [24] under the name of incremental-IOSS. □

## 2.1 iIOSS Lyapunov Functions

**Definition 2.5** We call a continuous function $V : \mathbb{R}^n \to \mathbb{R}_{\geq 0}$ an *iIOSS Lyapunov function* if there exist $\underline{\alpha}$, $\overline{\alpha} \in \mathcal{K}_\infty$, $\sigma_1$, $\sigma_2 \in \mathcal{K}$, and $\kappa : \mathbb{R}_{\geq 0} \to \mathbb{R}_{\geq 0}$ continuous positive definite so that

$$\underline{\alpha}(|\xi|) \leq V(\xi) \leq \overline{\alpha}(|\xi|) \qquad \forall \xi \in \mathbb{R}^n \tag{4}$$

and

$$\zeta \cdot f(\xi,\mu) \leq -\kappa(|\xi|) + \sigma_1(|h(\xi)|) + \sigma_2(|\mu|) \qquad \forall \xi \in \mathbb{R}^n,\ \forall \mu \in \mathbb{R}^m \tag{5}$$

for each $\zeta \in \partial_D V(\xi)$. □

We remark that the decrease statement (5) can be written equivalently in an integral formulation, using the following standard result.

**Proposition 2.6** (e.g. [20] Proposition 14) Given a forward complete system as in (1), a continuous function $V : \mathbb{R}^n \to \mathbb{R}_{\geq 0}$, and a continuous function $w : \mathbb{R}^n \times \mathbb{R}^m \to \mathbb{R}$, the following are equivalent:

1. For all $\xi \in \mathbb{R}^n$ and all $\mu \in \mathbb{R}^m$

$$\zeta \cdot f(\xi,\mu) \leq w(\xi,\mu)$$

for each $\zeta \in \partial_D V(\xi)$.

2. For each $\xi \in \mathbb{R}^n$ and each input $u$, the solution $x(\cdot,\xi,u)$ satisfies

$$V(x(t,\xi,u)) - V(\xi) \leq \int_0^t w(x(s,\xi,u),u(s))\,ds$$

for any $t \geq 0$.

□

**Remark 2.7** Applying Proposition 2.6 with

$$w(\xi, \mu) = -\kappa(|\xi|) + \sigma_1(|h(\xi)|) + \sigma_2(|\mu|),$$

we conclude that the decrease statement (5) in the definition of an iIOSS Lyapunov function could be equivalently written as

$$V(x(t, \xi, u)) - V(\xi) \leq \int_0^t -\kappa(|x(s, \xi, u)|) + \sigma_1(|h(x(s, \xi, u))|) + \sigma_2(|u(s)|)\, ds \tag{6}$$

for all $\xi \in \mathbb{R}^n$, all inputs $u$, and all $t \geq 0$. This alternative formulation will be used below. □

## 3 Lyapunov Characterization

Our main result is the following

**Theorem 1** *Suppose system (1) is forward complete. The following are equivalent.*

1. *The system is iIOSS.*
2. *The system admits an iIOSS Lyapunov function.*

**Remark 3.1** The main result in [11] is a Lyapunov characterization of the IOSS property. It is shown in that reference that a system is IOSS if and only if it admits an IOSS Lyapunov function, which can be defined as an iIOSS Lyapunov function for which the function $\kappa$ is of class $\mathcal{K}_\infty$. Thus it is an immediate consequence of Theorem 1 that the IOSS property implies the iIOSS property. □

**Remark 3.2** We will prove a slightly stronger statement than ($2 \Rightarrow 1$) of Theorem 1. The proof below shows that the existence of a *lower semicontinuous* iIOSS Lyapunov function implies that a system is iIOSS. □

It is still an open question whether every iIOSS system admits a *smooth* iIOSS Lyapunov function. However, as a minor extension of the proof of Theorem 1 we will also show the following.

**Lemma 3.3** Suppose system (1) is forward complete and has compact input value set $\mathbb{U}$. Then, if the system is iIOSS, it admits a *smooth* iIOSS Lyapunov function, i.e. there exists a smooth ($C^\infty$) function $V : \mathbb{R}^n \to \mathbb{R}_{\geq 0}$, and $\underline{\alpha}$, $\overline{\alpha} \in \mathcal{K}_\infty$, $\sigma_1$, $\sigma_2 \in \mathcal{K}$, and $\kappa : \mathbb{R}_{\geq 0} \to \mathbb{R}_{\geq 0}$ continuous positive definite so that (4) holds and

$$\nabla V(\xi) \cdot f(\xi, \mu) \leq -\kappa(|\xi|) + \sigma_1(|h(\xi)|) + \sigma_2(|\mu|) \qquad \forall \xi \in \mathbb{R}^n,\ \forall \mu \in \mathbb{U}.$$

**Remark 3.4** Lemma 3.3 provides, in particular, a Lyapunov characterization in terms of a smooth function for the property of *integral-output to state stability* (iOSS) which is defined as iIOSS for systems with no inputs (or equivalently, with $\mathbb{U} = \{0\}$). □

### 3.1 Sufficiency

We begin with the proof of (2 $\Rightarrow$ 1) (sufficiency) in Theorem 1. Here we follow the sufficiency argument given in [1]. A few preliminary lemmas are needed.

**Lemma 3.5** ([1] Lemma 4.1) Let $\kappa : \mathbb{R}_{\geq 0} \to \mathbb{R}_{\geq 0}$ be a continuous positive definite function. Then there exists $\rho_1 \in \mathcal{K}_\infty$ and $\rho_2 \in \mathcal{L}$ such that

$$\kappa(s) \geq \rho_1(s)\rho_2(s) \qquad \forall s \geq 0.$$

□

The following comparison result will be needed. This is a generalization of Corollary 4.3 in [1].

**Proposition 3.6** Given any continuous positive definite $\alpha : \mathbb{R}_{\geq 0} \to \mathbb{R}_{\geq 0}$, there exists a $\mathcal{KL}$ function $\beta$ with the following property. For any $0 < \widetilde{t} \leq \infty$, any lower semicontinuous function $y : [0, \widetilde{t}\,) \to \mathbb{R}_{\geq 0}$, and any measurable, locally essentially bounded function $v : [0, \widetilde{t}\,) \to \mathbb{R}_{\geq 0}$, if

$$y(t_2) \leq y(t_1) + \int_{t_1}^{t_2} -\alpha(y(s)) + v(s)\, ds \qquad \forall\, 0 \leq t_1 \leq t_2 < \widetilde{t}, \tag{7}$$

then

$$y(t) \leq \beta(y(0), t) + \int_0^t 2v(s)\, ds \qquad \forall t \in [0, \widetilde{t}\,).$$

□

The following lemma will be needed to prove Proposition 3.6.

**Lemma 3.7** Suppose given a locally Lipschitz positive definite function $\alpha : \mathbb{R}_{\geq 0} \to \mathbb{R}_{\geq 0}$, a time $0 < \widetilde{t} \leq \infty$, and a measurable, locally essentially bounded function $v : [0, \widetilde{t}\,) \to \mathbb{R}_{\geq 0}$. Let $y : [0, \widetilde{t}\,) \to \mathbb{R}_{\geq 0}$ be any lower semicontinuous function which satisfies (7). Define $w(\cdot)$ to be the solution of the initial value problem

$$\dot{w}(t) = -\alpha(w(t)) + v(t), \qquad w(0) = y(0). \tag{8}$$

Then $w(t)$ is defined for all $t \in [0, \widetilde{t}\,)$ and

$$y(t) \leq w(t) \qquad \forall t \in [0, \widetilde{t}\,).$$

*Proof.* (We follow the proof of Theorem III.4.1 in [5]). Let $y(\cdot)$ and $w(\cdot)$ be as above for given $\alpha$, $\widetilde{t}$, and $v(\cdot)$. We first note that $w(\cdot)$ exists for all $t \in [0, \widetilde{t}\,)$, since $\alpha$ is nonnegative and $v(\cdot)$ is essentially bounded on each finite interval. For each integer $n \geq 1$, let $w_n(\cdot)$ be the solution of

$$\dot{w}_n(t) = -\alpha(w_n(t)) + v(t) + \frac{1}{n}, \qquad w_n(0) = y(0), \tag{9}$$

which is also defined on $[0, \widetilde{t}\,)$. We will show that

$$y(t) \leq w_n(t) \qquad \forall t \in [0, \widetilde{t}\,) \tag{10}$$

for all $n \geq 1$. Indeed, suppose not. Then there exists $n \geq 1$ and $\tau \in [0, \widetilde{t}\,)$ so that

$$y(\tau) > w_n(\tau).$$

Let

$$t_0 := \sup \{t \in [0, \tau] \,:\, y(t) \leq w_n(t)\}.$$

Then, as $y(\cdot)$ is lower semicontinuous and $w_n(\cdot)$ is continuous,

$$y(t_0) \leq w_n(t_0).$$

We claim that in fact $y(t_0) = w_n(t_0)$. If this were not the case, there would be numbers $\delta_1$, $\delta_2$ so that

$$y(t_0) < \delta_1 < \delta_2 < w_n(t_0). \tag{11}$$

From (7), we have that

$$y(t_0 + t) \leq y(t_0) + \int_{t_0}^{t_0+t} -\alpha(y(s)) + v(s)\, ds$$

for each $t \in [0, \widetilde{t} - t_0)$. Since

$$\lim_{t \to 0} \int_{t_0}^{t_0+t} -\alpha(y(s)) + v(s)\, ds = 0,$$

it follows from (11) that there is some $\varepsilon_1 > 0$ so that $y(t_0 + t) < \delta_1$ for all $t \in [0, \varepsilon_1]$. Since $w_n(\cdot)$ is continuous, (11) also gives an $\varepsilon_2 > 0$ so that $w_n(t_0 + t) > \delta_2$ for all $t \in [0, \varepsilon_2]$. Thus $y(t_0 + t) < w_n(t_0 + t)$ for all $t$ sufficiently small, which contradicts the definition of $t_0$. We conclude that $y(t_0) = w_n(t_0)$.

From (7) and Taylor's Theorem, we have that for $\varepsilon \in [0, \widetilde{t} - t_0)$,

$$\begin{aligned} y(t_0 + \varepsilon) &\leq y(t_0) + \int_{t_0}^{t_0+\varepsilon} -\alpha(y(s)) + v(s)\, ds \\ &= y(t_0) - \varepsilon\alpha(y(t_0)) + \varepsilon v(t_0) + o(\varepsilon) \end{aligned}$$

and from (9)

$$\begin{aligned} w_n(t_0 + \varepsilon) &= w_n(t_0) + \int_{t_0}^{t_0+\varepsilon} -\alpha(w_n(s)) + v(s) + \frac{1}{n}\, ds \\ &= w_n(t_0) - \varepsilon\alpha(w_n(t_0)) + \varepsilon v(t_0) + \frac{\varepsilon}{n} + o(\varepsilon), \end{aligned}$$

where $o(\cdot)$ signifies a function satisfying $\lim_{t \to 0} \frac{o(t)}{t} = 0$. Since $w_n(t_0) = y(t_0)$, it follows that $y(t_0 + \varepsilon) \leq w_n(t_0 + \varepsilon)$ for $\varepsilon$ sufficiently small, a contradiction. Thus (10) holds for all $n \geq 1$.

We note that $w_n(t) \to w(t)$ uniformly on each finite time interval (cf. e.g. Theorem 1 in [25]). Thus for any $T \in [0, \widetilde{t}\,)$, as (10) holds for all $n$,

$$y(t) \leq \lim_{n \to \infty} w_n(t) = w(t) \qquad \forall t \in [0, T].$$

As $T > 0$ is arbitrary, we conclude that $y(t) \leq w(t)$ for all $t \in [0, \widetilde{t}\,)$. ∎

To complete the proof of Proposition 3.6 we will also need the following statement.

**Lemma 3.8** ([1] Corollary 4.3) Given any continuous positive definite $\alpha : \mathbb{R}_{\geq 0} \to \mathbb{R}_{\geq 0}$, there exists a $\mathcal{KL}$ function $\beta$ with the following property. For any $0 < \widetilde{t} \leq \infty$, any absolutely continuous function $w : [0, \widetilde{t}\,) \to \mathbb{R}_{\geq 0}$, and any measurable, locally essentially bounded function $v : [0, \widetilde{t}\,) \to \mathbb{R}_{\geq 0}$, if

$$\dot{w}(t) \leq -\alpha(w(t)) + v(t) \tag{12}$$

for almost all $t \in [0, \widetilde{t}\,)$, then

$$w(t) \leq \beta(w(0), t) + \int_0^t 2v(s)\, ds \qquad \forall t \in [0, \widetilde{t}\,).$$

□

The proof of Proposition 3.6 is a straightforward combination of Lemma 3.7 and Lemma 3.8.
*Proof.* (Proposition 3.6) Let a continuous positive definite $\alpha : \mathbb{R}_{\geq 0} \to \mathbb{R}_{\geq 0}$ be given. Without loss of generality, we assume $\alpha$ is locally Lipschitz (otherwise we replace $\alpha$ by a locally Lipschitz function majorized by $\alpha$). Let $\beta$ be the $\mathcal{KL}$ function given by Lemma 3.8. Suppose $\widetilde{t}$, $y(\cdot)$ and $v(\cdot)$ are as in the statement of the Proposition so that (7) holds. Let $w(\cdot)$ be the solution of the initial value problem (8). Then Lemma 3.7 gives

$$y(t) \leq w(t) \qquad \forall t \in [0, \widetilde{t}\,).$$

Also, since $w(\cdot)$ satisfies (12) (as an equality), Lemma 3.8 gives

$$w(t) \leq \beta(w(0), t) + \int_0^t 2v(s)\, ds \qquad \forall t \in [0, \widetilde{t}\,).$$

Since $w(0) = y(0)$, the result follows. ■

We can now give the argument for sufficiency of the Lyapunov characterization. As mentioned earlier, this proof holds for lower semicontinuous Lyapunov functions.
*Proof.* Theorem 1 ($2 \Rightarrow 1$)
Suppose the function $V$ satisfies the definition of an iIOSS Lyapunov function for the forward complete system (1) with functions $\underline{\alpha}$, $\overline{\alpha}$, $\kappa$, $\sigma_1$ and $\sigma_2$ satisfying (4) and (5). Let $\rho_1 \in \mathcal{K}_\infty$ and $\rho_2 \in \mathcal{L}$ be functions as in Lemma 3.5 for $\kappa$. Let

$$\widetilde{\rho}(s) := \rho_1(\overline{\alpha}^{-1}(s))\rho_2(\underline{\alpha}^{-1}(s)).$$

By (4) and (6), we have, for each $\xi \in \mathbb{R}^n$ and each input $u$,

$$\begin{aligned}
V(x(t_2, \xi, u)) \;&\leq\; V(x(t_1, \xi, u)) + \int_{t_1}^{t_2} -\kappa(|x(s, \xi, u)|) + \sigma_1(|h(x(s, \xi, u))|) + \sigma_2(|u(s)|)\, ds \\
&\leq\; V(x(t_1, \xi, u)) + \int_{t_1}^{t_2} -\rho_1(|x(s, \xi, u)|)\rho_2(|x(s, \xi, u)|) \\
&\qquad\qquad + \sigma_1(|h(x(s, \xi, u))|) + \sigma_2(|u(s)|)\, ds \\
&\leq\; V(x(t_1, \xi, u)) + \int_{t_1}^{t_2} -\widetilde{\rho}(V(x(s, \xi, u))) + \sigma_1(|h(x(s, \xi, u))|) + \sigma_2(|u(s)|)\, ds
\end{aligned}$$

for all $0 \leq t_1 \leq t_2$.

Then, as $\widetilde{\rho}$ is continuous positive definite, Proposition 3.6 gives the existence of a $\mathcal{KL}$ function $\beta$ so that for each $\xi \in \mathbb{R}^n$ and each input $u$

$$\begin{aligned}\underline{\alpha}(|x(t,\xi,u)|) \;\le\; V(x(t,\xi,u)) \;&\le\; \beta(V(\xi),t) + \int_0^t 2\sigma_1(|h(x(s,\xi,u))|) + 2\sigma_2(|u(s)|)\,ds \\ &\le\; \beta(\overline{\alpha}(|\xi|),t) + \int_0^t 2\sigma_1(|h(x(s,\xi,u))|) + 2\sigma_2(|u(s)|)\,ds\end{aligned}$$

for all $t \ge 0$, which is the required bound. ∎

## 3.2 Necessity

We next prove (1 ⇒ 2) (necessity) for Theorem 1. We will construct an iIOSS Lyapunov function for a given iIOSS system. The proof combines ideas from the constructions in [28] and [1]. The following result will be needed.

This statement follows directly from Proposition 7 in [26].

**Proposition 3.9** For any given $\mathcal{KL}$ function $\beta$, there exist a family of mappings $\{T_r\}_{r\ge 0}$ with:

- for each fixed $r > 0$, $T_r : \mathbb{R}_{>0} \overset{\text{onto}}{\to} \mathbb{R}_{>0}$ is strictly decreasing;
- for each fixed $\varepsilon > 0$, $T_r(\varepsilon)$ is strictly increasing as $r$ increases and $\lim_{r\to\infty} T_r(\varepsilon) = \infty$;
- the map $(r,\varepsilon) \mapsto T_r(\varepsilon)$ is jointly continuous in $r$ and $\varepsilon$;

such that

$$\beta(s,t) \le \varepsilon$$

for all $s \le r$, all $t \ge T_r(\varepsilon)$. □

Before giving the construction, we will cite a lemma on boundedness of reachable sets for forward complete systems which says that the reachable set from a given point over a finite time interval $[0,T]$ is bounded if the inputs are required to satisfy a bound of the type

$$\int_0^T \gamma(|u(s)|)\,ds \le M < \infty \tag{13}$$

for an appropriate choice of $\gamma \in \mathcal{K}_\infty$.

**Remark 3.10** Note that for arbitrary $\mathcal{K}_\infty$ functions $\gamma$, this need not hold. Take, for example the one-dimensional system $\dot{x} = u^2$. With $\gamma(s) = s$, the inputs

$$u_k(t) = \begin{cases} k & 0 \le t \le \frac{1}{k} \\ 0 & \frac{1}{k} < t \le 1 \end{cases}$$

defined on $[0,1]$ satisfy $\int_0^1 \gamma(|u_k(s)|)\,ds = 1$ for each $k \ge 1$. However, the solution starting at the origin corresponding to the input $u_k(\cdot)$ satisfies $x(1) = k$, and so clearly one can reach an unbounded set in one time unit using controls satisfying (13). □

The following lemma shows that one can always choose $\gamma$ so that the bound (13) on inputs implies a bounded reachable set. (In the example above, clearly $\gamma(s) = s^2$ will do.)

**Lemma 3.11** ([2] Corollary 2.13) Suppose system (1) is forward complete. Then there exist functions $\chi_1$, $\chi_2$, $\chi_3$, $\sigma$ of class $\mathcal{K}_\infty$ and a constant $c \geq 0$ such that

$$|x(t,\xi,u)| \leq \chi_1(t) + \chi_2(|\xi|) + \chi_3\left(\int_0^t \sigma(|u(s)|)\,ds\right) + c$$

holds for all $\xi \in \mathbb{R}^n$, all inputs $u$, and all $t \geq 0$. □

We now provide the Lyapunov construction.

*Proof.* Theorem 1 (1 ⇒ 2)

Suppose the system (1) is forward complete and satisfies the iIOSS property with gains $\alpha$, $\beta$, $\gamma_1$ and $\gamma_2$.

Pick any smooth, strictly increasing and bounded function $k : \mathbb{R} \to \mathbb{R}_{>0}$ whose derivative is strictly decreasing. Then there are two positive numbers $c_1 < c_2$ so that $k(t) \in [c_1, c_2]$ for all $t \geq 0$. Define $\lambda(t) : \mathbb{R}_{\geq 0} \to \mathbb{R}_{>0}$ by

$$\lambda(t) := \frac{d}{dt}k(t).$$

Since the system is forward complete, we may find a function $\sigma \in \mathcal{K}_\infty$ as in Lemma 3.11. Define $\widetilde{\gamma}_2(s) := \max\{\gamma_2(s), \sigma(s)\}$ for all $s \geq 0$. Note that the iIOSS bound (2) holds with $\widetilde{\gamma}_2$ in the place of $\gamma_2$.

We define a Lyapunov function as

$$V_0(\xi) := \sup_u \sup_{t\geq 0}\left\{\left(\alpha(|x(t,\xi,u)|) - \int_0^t \gamma_1(|y(s,\xi,u)|)\,ds - \int_0^\infty 2\widetilde{\gamma}_2(|u(s)|)\,ds\right)k(t)\right\}$$

for each $\xi \in \mathbb{R}^n$. It is immediate that this function satisfies (4), as

$$c_1\alpha(|\xi|) \leq V_0(\xi) \leq c_2\beta(|\xi|, 0) \qquad \forall \xi \in \mathbb{R}^n. \tag{14}$$

The first of these inequalities follows from considering the trajectory with input $u \equiv 0$ at time $t = 0$, and the second from the iIOSS bound (3): for any $\xi \in \mathbb{R}^n$ and any input $u$,

$$\begin{aligned}
&\alpha(|x(t,\xi,u)|) - \int_0^t \gamma_1(|y(s,\xi,u)|)\,ds - \int_0^\infty 2\widetilde{\gamma}_2(|u(s)|)\,ds \\
\leq\; & \alpha(|x(t,\xi,u)|) - \int_0^t \gamma_1(|y(s,\xi,u)|)\,ds - \int_0^\infty \gamma_2(|u(s)|)\,ds \\
\leq\; & \beta(|\xi|, t) \\
\leq\; & \beta(|\xi|, 0)
\end{aligned} \tag{15}$$

for all $t \geq 0$.

Next, we observe that for each $\xi$, the supremum over inputs in $V_0(\xi)$ can be taken to be a supremum over a restricted set, as follows. From the iIOSS bound (3), we have, for any $\xi \in \mathbb{R}^n$ and any input $u$,

$$\alpha(|x(t,\xi,u)|) \leq \beta(|\xi|, 0) + \int_0^t \gamma_1(|y(s,\xi,u)|)\,ds + \int_0^\infty \widetilde{\gamma}_2(|u(s)|)\,ds$$

for all $t \geq 0$. Suppose now that $\xi$ and $u$ are such that

$$\int_0^\infty \widetilde{\gamma}_2(|u(s)|)\,ds > \beta(|\xi|, 0).$$

In this case it follows that

$$\alpha(|x(t,\xi,u)|) \le \int_0^t \gamma_1(|y(s,\xi,u)|)\,ds + \int_0^\infty 2\widetilde{\gamma}_2(|u(s)|)\,ds$$

for all $t \ge 0$. Then for all $t \ge 0$

$$\alpha(|x(t,\xi,u)|) - \int_0^t \gamma_1(|y(s,\xi,u)|)\,ds - \int_0^\infty 2\widetilde{\gamma}_2(|u(s)|)\,ds \le 0.$$

Since $V_0(\xi) \ge 0$ for each $\xi \in \mathbb{R}^n$, it follows that for each $\xi \in \mathbb{R}^n$

$$V_0(\xi) = \sup_{u\in U(|\xi|)} \sup_{t\ge 0} \left\{ \left( \alpha(|x(t,\xi,u)|) - \int_0^t \gamma_1(|y(s,\xi,u)|)\,ds - \int_0^\infty 2\widetilde{\gamma}_2(|u(s)|)\,ds \right) k(t) \right\}.$$

where, for each $r \ge 0$, we define $U(r) := \{u(\cdot)\ :\ \int_0^\infty \widetilde{\gamma}_2(|u(s)|)\,ds \le \beta(r,0)\}$.

We next make the observation that the supremum in time can be taken over a restricted set as well. Let $T_r(\varepsilon)$ be defined as in Proposition 3.9 for the function $\beta$. From (14) and (15) we have

$$V_0(\xi) = \sup_{u\in U(|\xi|)} \sup_{0\le t\le T_\xi} \left\{ \left( \alpha(|x(t,\xi,u)|) - \int_0^t \gamma_1(|y(s,\xi,u)|)\,ds - \int_0^t 2\widetilde{\gamma}_2(|u(s)|)\,ds \right) k(t) \right\},$$

where for each $\xi \in \mathbb{R}^n$ we set $T_\xi := T_{2|\xi|}(\frac{c_1}{c_2}\alpha(\frac{|\xi|}{2}))$.

We will show that the function $V_0$ is continuous on $\mathbb{R}^n$ by showing lower and upper semicontinuity in the next two lemmas.

**Proposition 3.12** The function $V_0$ is lower semicontinuous on $\mathbb{R}^n$.

*Proof.* We will show

$$\liminf_{\xi\to\xi_0} V_0(\xi) \ge V_0(\xi_0)$$

for all $\xi_0 \in \mathbb{R}^n$.

Fix $\xi_0 \in \mathbb{R}^n$ and let $\varepsilon > 0$ be given. There exists an input $u_0$ and a time $t_0 \ge 0$ so that

$$\left( \alpha(|x(t_0,\xi_0,u_0)|) - \int_0^{t_0} \gamma_1(|y(s,\xi_0,u_0)|)\,ds - \int_0^\infty 2\widetilde{\gamma}_2(|u_0(s)|)\ ds \right) k(t_0) \ge V_0(\xi_0) - \frac{\varepsilon}{2}.$$

By continuity of $x(t_0,\cdot,u_0)$ and $\alpha$, there exists a neighbourhood $U_1$ of $\xi_0$ so that

$$|\alpha(|x(t_0,\xi_0,u_0)|) - \alpha(|x(t_0,\xi,u_0)|)| \le \frac{\varepsilon}{4k(t_0)}$$

for all $\xi \in U_1$. Furthermore, as $\xi \to \xi_0$ implies $y(t,\xi,u_0)$ converges uniformly to $y(t,\xi_0,u_0)$ on the finite interval $[0,t_0]$, and since $\gamma_1$ is uniformly continuous on a compact containing an open neighbourhood of $\{y(t,\xi_0,u_0)\ :\ t\in[0,t_0]\}$, we can find a neighbourhood $U_2 \subseteq U_1$ of $\xi_0$ so that each $\xi \in U_2$ satisfies

$$|\gamma_1(|y(s,\xi_0,u_0)|) - \gamma_1(|y(s,\xi,u_0)|)| \le \frac{\varepsilon}{4t_0k(t_0)}$$

for all $s \in [0, t_0]$. Then for each $\xi \in U_2$,

$$
\begin{aligned}
& \left| \left( \alpha(|x(t_0, \xi_0, u_0)|) - \int_0^{t_0} \gamma_1(|y(s, \xi_0, u_0)|)\, ds \right) - \left( \alpha(|x(t_0, \xi, u_0)|) - \int_0^{t_0} \gamma_1(|y(s, \xi, u_0)|)\, ds \right) \right| \\
\leq \quad & \frac{\varepsilon}{4k(t_0)} + \int_0^{t_0} |\gamma_1(|y(s, \xi_0, u_0)|) - \gamma_1(|y(s, \xi, u_0)|)|\, ds \\
\leq \quad & \frac{\varepsilon}{2k(t_0)}.
\end{aligned}
$$

This gives, for each $\xi \in U_2$,

$$
\begin{aligned}
V_0(\xi) \quad \geq \quad & \left( \alpha(|x(t_0, \xi, u_0)|) - \int_0^{t_0} \gamma_1(|y(s, \xi, u_0)|)\, ds - \int_0^{\infty} 2\widetilde{\gamma}_2(|u_0(s)|)\, ds \right) k(t_0) \\
\geq \quad & \left( \alpha(|x(t_0, \xi_0, u_0)|) - \int_0^{t_0} \gamma_1(|y(s, \xi_0, u_0)|)\, ds - \int_0^{\infty} 2\widetilde{\gamma}_2(|u_0(s)|)\, ds \right) k(t_0) - \frac{\varepsilon}{2} \\
\geq \quad & V_0(\xi_0) - \varepsilon,
\end{aligned}
$$

Hence $V_0$ is lower semicontinuous. ∎

The next result will be needed to show upper semicontinuity.

**Proposition 3.13** For each $T > 0$ and each compact $C \subset \mathbb{R}^n$, there exists $L_{C,T} > 0$ so that for any input $u \in U(C) := \bigcup_{\xi \in C} U(|\xi|)$, each pair $\eta$, $\zeta \in C$ has the property that

$$
|x(t, \eta, u) - x(t, \zeta, u)| \leq L_{C,T}\, |\eta - \zeta| \qquad \forall t \in [0, T].
$$

That is, the trajectories are Lipschitz in the initial conditions, uniformly over inputs $u \in U(C)$.

*Proof.* From Lemma 3.11, we have that the trajectories stay in a bounded set on the interval $[0, T]$. A standard Gronwall's Lemma argument gives this Lipschitz condition from the local Lipschitz assumption on $f$. ∎

**Proposition 3.14** The function $V_0$ is upper semicontinuous on $\mathbb{R}^n$.

*Proof.* We will show

$$
\limsup_{\xi \to \xi_0} V_0(\xi) \leq V_0(\xi_0) \tag{16}
$$

for all $\xi_0 \in \mathbb{R}^n$.

Suppose (16) fails at some $\xi_0 \in \mathbb{R}^n$. Then there exists $\varepsilon > 0$ and a sequence $\{\xi_j\}_{j=1}^{\infty}$ so that $\xi_j \to \xi_0$ and

$$
V_0(\xi_j) > V_0(\xi_0) + \varepsilon \tag{17}
$$

for all $j \geq 1$. Choose $r > 0$ so that $|\xi_0| \leq r$ and $|\xi_j| \leq r$ for all $j \geq 1$. Then for each $j \geq 1$,

$$
V(\xi_j) = \sup_{u \in U(r)} \max_{t \in [0, T_0]} \left\{ \left( \alpha(|x(t, \xi_j, u)|) - \int_0^{t} \gamma_1(|y(s, \xi_j, u)|)\, ds - \int_0^{\infty} 2\widetilde{\gamma}_2(|u(s)|)\, ds \right) k(t) \right\}
$$

where $T_0 := T_r(\frac{1}{c_2}(V_0(\xi_0)+\varepsilon))$. Now, for each $j \geq 1$, there exists $\tau_j \in [0, T_0]$, $u_j \in U(r)$ for which

$$\left(\alpha(|x(\tau_j,\xi_j,u_j)|) - \int_0^{\tau_j} \gamma_1(|y(s,\xi_j,u_j)|)\,ds - \int_0^{\infty} 2\widetilde{\gamma}_2(|u_j(s)|)\,ds\right) k(\tau_j) \geq V_0(\xi_j) - \frac{\varepsilon}{2}.$$

Let $R$ be the reachable set from $B_r := \{\xi \in \mathbb{R}^n : |\xi| \leq r\}$ in time less than or equal to $T_0$ with controls in $U(r)$. Then Lemma 3.11 tells us that $R$ is bounded. From Proposition 3.13 and the fact that the output map $h$ is locally Lipschitz, we can find $L_x > 0$ and $L_y > 0$ so that

$$\begin{aligned} |x(t,\eta,u) - x(t,\zeta,u)| &\leq L_x\,|\eta - \zeta| \\ |y(t,\eta,u) - y(t,\zeta,u)| &\leq L_y\,|\eta - \zeta| \end{aligned}$$

for all $t \in [0, T_0]$, for any pair $\eta$, $\zeta \in B_r$ and for any input $u \in U(r)$. Further, since $\alpha$ is continuous, it is uniformly continuous on the bounded set $R$, so there is a $\delta_x > 0$ so that

$$|\zeta_1 - \zeta_2| \leq \delta_x \;\Rightarrow |\alpha(\zeta_1) - \alpha(\zeta_2)| \leq \frac{\varepsilon}{4k(T_0)}$$

if $\zeta_1$, $\zeta_2 \in R$. Likewise, since $\gamma_1$ is uniformly continuous on the bounded set $h(R) := \{h(\eta) : \eta \in R\}$, there is a $\delta_y > 0$ so that

$$|\eta_1 - \eta_2| \leq \delta_y \;\Rightarrow |\gamma_1(\eta_1) - \gamma_1(\eta_2)| \leq \frac{\varepsilon}{4T_0 k(T_0)}$$

for $\eta_1$, $\eta_2 \in h(R)$. Then, for each $j$ large enough so that $|\xi_j - \xi_0| \leq \frac{\min\{\delta_x,\delta_y\}}{\max\{L_x,L_y\}}$, we find

$$\begin{aligned} &\left|\left(\alpha(|x(\tau_j,\xi_j,u_j)|) - \int_0^{\tau_j} \gamma_1(|y(s,\xi_j,u_j)|)\,ds\right) - \left(\alpha(|x(\tau_j,\xi_0,u_j)|) - \int_0^{\tau_j} \gamma_1(|y(s,\xi_0,u_j)|)\,ds\right)\right| \\ &\leq \frac{\varepsilon}{4k(T_0)} + \int_0^{\tau_j} |\gamma_1(|y(s,\xi_j,u_j)|) - \gamma_1(|y(s,\xi_0,u_j)|)|\,ds \\ &\leq \frac{\varepsilon}{2k(T_0)}. \end{aligned}$$

From which we find, for each $j$ sufficiently large,

$$\begin{aligned} V_0(\xi_0) &\geq \left(\alpha(|x(\tau_j,\xi_0,u_j)|) - \int_0^{\tau_j} \gamma_1(|y(s,\xi_0,u_j)|)\,ds - \int_0^{\infty} 2\widetilde{\gamma}_2(|u_j(s)|)\,ds\right) k(\tau_j) \\ &\geq \left(\alpha(|x(\tau_j,\xi_j,u_j)|) - \int_0^{\tau_j} \gamma_1(|y(s,\xi_j,u_j)|)\,ds - \int_0^{\infty} 2\widetilde{\gamma}_2(|u_j(s)|)\,ds\right) k(\tau_j) - \frac{\varepsilon}{2} \\ &\geq V_0(\xi_j) - \varepsilon \end{aligned}$$

which contradicts (17). We conclude that $V_0$ is upper semicontinuous. ∎

Finally, we show that the function $V_0$ satisfies the decrease statement (5). Let $\xi \in \mathbb{R}^n \backslash \{0\}$ and an input $v$ be given, and consider the resulting trajectory. For $\tau > 0$ small enough, we have $\frac{|\xi|}{2} < |x(\tau,\xi,v)| < 2\,|\xi|$, so for such $\tau$ the supremum over time in the expression for $V_0(x(\tau,\xi,v))$ may be taken over $[0, T_\xi]$. We find, for such $\tau$ sufficiently small,

$$\begin{aligned} &V_0(x(\tau,\xi,v)) \\ &= \sup_u \sup_{0 \leq s \leq T_\xi} \left\{\left(\alpha(|x(s,x(\tau,\xi,v),u)|) - \int_0^{s} \gamma_1(|y(r,x(\tau,\xi,v),u)|)\,dr\right.\right. \end{aligned}$$

$$
\begin{aligned}
& \qquad\qquad\qquad\qquad\qquad\qquad - \int_0^\infty 2\widetilde{\gamma}_2(|u(r)|)\,dr\Big)\,k(s)\Big\} \\
= & \sup_u \sup_{\tau\le t\le \tau+T_\xi} \Big\{\Big(\alpha(|x(t,\xi,v\sharp_\tau u)|) - \int_\tau^t \gamma_1(|y(r,\xi,v\sharp_\tau u)|)\,dr \\
& \qquad\qquad\qquad\qquad - \int_\tau^\infty 2\widetilde{\gamma}_2(|v\sharp_\tau u(r)|)\,dr\Big)\,k(t-\tau)\Big\} \\
\le & \sup_u \sup_{0\le t\le \tau+T_\xi} \Big\{\Big(\alpha(|x(t,\xi,v\sharp_\tau u)|) - \int_0^t \gamma_1(|y(r,\xi,v\sharp_\tau u)|)\,dr + \int_0^\tau \gamma_1(|y(r,\xi,v\sharp_\tau u)|)\,dr \\
& \qquad - \int_0^\infty 2\widetilde{\gamma}_2(|v\sharp_\tau u(r)|)\,dr + \int_0^\tau 2\widetilde{\gamma}_2(|v\sharp_\tau u(r)|)\,dr\Big)\,k(t-\tau)\Big\} \\
\le & \sup_u \sup_{0\le t\le \tau+T_\xi} \Big\{\Big(\alpha(|x(t,\xi,v\sharp_\tau u)|) - \int_0^t \gamma_1(|y(r,\xi,v\sharp_\tau u)|)\,dr \\
& \qquad - \int_0^\infty 2\widetilde{\gamma}_2(|v\sharp_\tau u(r)|)\,dr\Big)\,k(t-\tau)\Big\} + c_2\int_0^\tau \gamma_1(|y(r,\xi,v)|) + 2\widetilde{\gamma}_2(|v(r)|)\ dr \\
\le & \sup_u \sup_{0\le t\le \tau+T_\xi} \Big\{\Big(\alpha(|x(t,\xi,v\sharp_\tau u)|) - \int_0^t \gamma_1(|y(r,\xi,v\sharp_\tau u)|)\,dr \\
& \qquad - \int_0^\infty 2\widetilde{\gamma}_2(|v\sharp_\tau u(r)|)\,dr\Big)\,k(t)\Big\}\cdot \max_{0\le t\le \tau+T_\xi}\left\{\frac{k(t-\tau)}{k(t)}\right\} \\
& \qquad + c_2\int_0^\tau \gamma_1(|y(r,\xi,v)|) + 2\widetilde{\gamma}_2(|v(r)|)\ dr \\
\le & \sup_{\widehat{u}} \sup_{0\le t} \Big\{\Big(\alpha(|x(t,\xi,\widehat{u})|) - \int_0^t \gamma_1(|y(r,\xi,\widehat{u})|)\,dr - \int_0^\infty 2\widetilde{\gamma}_2(|\widehat{u}(r)|)\,dr\Big)\,k(t)\Big\}\cdot \\
& \qquad \max_{0\le t\le \tau+T_\xi}\left\{\frac{k(t-\tau)}{k(t)}\right\} + c_2\int_0^\tau \gamma_1(|y(r,\xi,v)|) + 2\widetilde{\gamma}_2(|v(r)|)\ dr \\
= & V_0(\xi)\cdot \max_{0\le t\le \tau+T_\xi}\left\{\frac{k(t-\tau)}{k(t)}\right\} + c_2\int_0^\tau \gamma_1(|y(r,\xi,v)|) + 2\widetilde{\gamma}_2(|v(r)|)\ dr,
\end{aligned}
$$

where $v\sharp_\tau u$ is the concatenation of $u$ with $v$ at time $\tau$, that is

$$
v\sharp_\tau u = \begin{cases} v(t) & \text{if } \ 0\le t\le \tau \\ u(t-\tau) & \text{if } \ \tau < t \end{cases} .
$$

Rewriting, we arrive at

$$
\begin{aligned}
& V_0(x(\tau,\xi,v)) - V_0(\xi) \\
\le \ & V_0(\xi)\cdot \max_{0\le t\le \tau+T_\xi}\left[-1 + \frac{k(t-\tau)}{k(t)}\right] + c_2\int_0^\tau \gamma_1(|y(r,\xi,v)|) + 2\widetilde{\gamma}_2(|v(r)|)\,dr
\end{aligned}
$$

for $\tau > 0$ sufficiently small. Recall that $\lambda(t) = \frac{d}{dt}k(t)$ is decreasing, so for $\tau > 0$ small enough,

$$
\begin{aligned}
\max_{0\le t\le \tau+T_\xi}\left[-1 + \frac{k(t-\tau)}{k(t)}\right] & \le -\min_{0\le t\le \tau+T_\xi}\left[\frac{k(t)-k(t-\tau)}{c_2}\right] \\
& = -\frac{k(T_\xi+\tau) - k(T_\xi)}{c_2}.
\end{aligned}
$$

So for $\tau > 0$ sufficiently small,

$$\begin{aligned} V_0(x(\tau,\xi,v)) - V_0(\xi) &\leq -V_0(\xi)\frac{k(T_\xi+\tau)-k(T_\xi)}{c_2} + c_2\int_0^\tau \gamma_1(|y(r,\xi,v)|) + 2\widetilde{\gamma}_2(|v(r)|)\,dr \\ &= \int_0^\tau -\frac{V_0(\xi)}{c_2}\lambda(T_\xi+r) + c_2[\gamma_1(|y(r,\xi,v)|) + 2\widetilde{\gamma}_2(|v(r)|)]\,dr. \qquad (18) \end{aligned}$$

Recall that (18) has been verified for all $\xi \neq 0$. We next note that it also holds for $\xi = 0$. Let an input $v$ be given. The calculation above (with $T_\xi = \infty$) gives, for $\tau > 0$,

$$\begin{aligned} V_0(x(\tau,0,v)) &\leq V_0(0)\cdot\sup_{0\leq t}\left\{\frac{k(t-\tau)}{k(t)}\right\} + c_2\int_0^\tau \gamma_1(|y(r,\xi,v)|) + 2\widetilde{\gamma}_2(|v(r)|)\;dr \\ &= c_2\int_0^\tau \gamma_1(|y(r,\xi,v)|) + 2\widetilde{\gamma}_2(|v(r)|)\;dr, \end{aligned}$$

as $V_0(0) = 0$. Clearly this gives (18) for $\xi = 0$.

Finally, we will make use of the following lemma to formulate the decrease statement (18) in the viscosity sense.

**Lemma 3.15** Suppose given a system as in (1), a function $V : \mathbb{R}^n \to \mathbb{R}$, a point $\xi \in \mathbb{R}^n$, and an element $\mu \in \mathbb{R}^m$. Then, if there exists a continuous $\alpha_{\xi,\mathbf{u}} : \mathbb{R}_{\geq 0} \to \mathbb{R}$ and $\varepsilon > 0$ so that for all $0 \leq \tau < \varepsilon$

$$V(x(\tau,\xi,\mathbf{u})) - V(\xi) \leq \int_0^\tau \alpha_{\xi,\mathbf{u}}(r)\,dr \qquad (19)$$

where $\mathbf{u}$ is the input constantly equal to $\mu$, then for any $\zeta \in \partial_D V(\xi)$, the instantaneous form of (19) holds in the viscosity sense:

$$\zeta\cdot f(\xi,\mu) \leq \alpha_{\xi,\mathbf{u}}(0).$$

*Proof.* Suppose $V$, $\xi$, $\mu$, $\alpha_{\xi,\mathbf{u}}$ and $\varepsilon$ are as above, and suppose $\zeta \in \partial_D V(\xi)$. Then, from the definition of the viscosity subgradient, we know that for $\tau$ small enough

$$\int_0^\tau \alpha_{\xi,\mathbf{u}}(r)\,dr \geq V(x(\tau,\xi,\mathbf{u})) - V(\xi) \geq \zeta\cdot(x(\tau,\xi,\mathbf{u})-\xi) + g(x(\tau,\xi,\mathbf{u})-\xi) \qquad (20)$$

where $g$ is some function satisfying $\lim_{s\to 0}\frac{g(s)}{|s|} = 0$. We note that since $\mathbf{u}$ is constant valued, the trajectory $x(\cdot,\xi,\mathbf{u})$ is differentiable (not merely absolutely continuous). In particular, $\frac{d}{dt}x(t,\xi,\mathbf{u})|_{t=0} = f(\xi,\mu)$. Now, dividing by $\tau$ in (20) gives

$$\frac{\int_0^\tau \alpha_{\xi,\mathbf{u}}(r)\,dr}{\tau} \geq \zeta\cdot\frac{x(\tau,\xi,\mathbf{u})-\xi}{\tau} + \frac{g(x(\tau,\xi,\mathbf{u})-\xi)}{\tau}.$$

Taking the limit as $\tau$ tends to 0, we find

$$\begin{aligned} \alpha_{\xi,\mathbf{u}}(0) &\geq \zeta\cdot f(\xi,\mu) + \lim_{\tau\to 0}\frac{g(x(\tau,\xi,\mathbf{u})-\xi)}{\tau} \\ &= \zeta\cdot f(\xi,\mu) + \lim_{\tau\to 0}\frac{|x(\tau,\xi,\mathbf{u})-\xi|}{\tau}\,\frac{g(x(\tau,\xi,\mathbf{u})-\xi)}{|x(\tau,\xi,\mathbf{u})-\xi|} \\ &= \zeta\cdot f(\xi,\mu) + |f(\xi,\mu)|\,\lim_{\tau\to 0}\frac{g(x(\tau,\xi,\mathbf{u})-\xi)}{|x(\tau,\xi,\mathbf{u})-\xi|} \\ &= \zeta\cdot f(\xi,\mu). \end{aligned}$$

■

Applying Lemma 3.15 to (18) with

$$\alpha_{\xi,\mathbf{u}}(r) = -\frac{V_0(\xi)}{c_2}\lambda(T_\xi + r) + c_2[\gamma_1(|y(r,\xi,\mathbf{u})|) + 2\widetilde{\gamma}_2(|\mathbf{u}(r)|)]$$

we conclude that

$$\zeta \cdot f(\xi,\mu) \leq -\frac{V_0(\xi)\lambda(T_\xi)}{c_2} + c_2\gamma_1(|h(\xi)|) + 2c_2\widetilde{\gamma}_2(|\mu|) \quad \forall \xi \in \mathbb{R}^n,\ \forall \mu \in \mathbb{R}^m$$

for each $\zeta \in \partial_D V_0(\xi)$. This implies (5), since $V_0(\cdot)\lambda(T_{(\cdot)}) : \mathbb{R}^n \to \mathbb{R}_{\geq 0}$ is continuous positive definite, so we can choose a continuous positive definite function $\kappa : \mathbb{R}_{\geq 0} \to \mathbb{R}_{\geq 0}$ so that $\kappa(|\xi|) \leq \frac{1}{c_2}V_0(\xi)\lambda(T_\xi)$ for each $\xi \in \mathbb{R}^n$. Then

$$\zeta \cdot f(\xi,\mu) \leq -\kappa(|\xi|) + c_2\gamma_1(|h(\xi)|) + 2c_2\widetilde{\gamma}_2(|\mu|) \quad \forall \xi \in \mathbb{R}^n,\ \mu \in \mathbb{R}^m$$

for each $\zeta \in \partial_D V_0(\xi)$.

This completes the construction of the iIOSS Lyapunov function. ■

Finally, we prove Lemma 3.3 by extending the proof in the case where the input value set $\mathbb{U}$ is compact.

*Proof.* (Lemma 3.3)

Suppose the forward complete system (1) satisfies the iIOSS property and has compact input value set $\mathbb{U}$. Given the construction above, it follows from Corollary 4.22 in [11] that there exists a smooth function $\widetilde{V} : \mathbb{R}^n \backslash \{0\} \to \mathbb{R}_{\geq 0}$, for which

$$\frac{V_0(\xi)}{2} \leq \widetilde{V}(\xi) \leq 2V_0(\xi) \qquad \forall \xi \in \mathbb{R}^n \backslash \{0\} \tag{21}$$

and

$$\nabla\widetilde{V}(\xi) \cdot f(\xi,\mu) \leq -\frac{1}{2}\kappa(|\xi|) + c_2\gamma_1(|h(\xi)|) + 2c_2\widetilde{\gamma}_2(|\mu|) \quad \forall \xi \in \mathbb{R}^n \backslash \{0\}, \forall \mu \in \mathbb{U}. \tag{22}$$

We extend $\widetilde{V}$ to $\mathbb{R}^n$ by setting $\widetilde{V}(0) = 0$. It is immediate that the resulting function is continuous on $\mathbb{R}^n$ and that (21) holds on all of $\mathbb{R}^n$.

By Proposition 4.2 in [15] there is a smooth $\rho \in \mathcal{K}_\infty$ with $\rho'(s) > 0$ for all $s > 0$ such that $\rho \circ \widetilde{V}$ is smooth everywhere. Without loss of generality, we may assume that $\rho'(s) \leq 1$ for all $s > 0$. (If it is not, we may replace $\rho$ by a smooth $\mathcal{K}_\infty$ function $\rho_0$ with the property that $\rho_0'(s) = \rho'(s)$ in a neighbourhood of the origin where $\rho'(s) \leq 1$ and $\rho_0'(s) \leq 1$ everywhere else.) Let $V = \rho \circ \widetilde{V}$. It follows from (21) and (14) that

$$\underline{\alpha}(|\xi|) \leq V(\xi) \leq \overline{\alpha}(|\xi|) \qquad \forall \xi \in \mathbb{R}^n,$$

where $\underline{\alpha}(s) = \rho(\frac{c_1}{2}\alpha(s))$ and $\overline{\alpha}(s) = \rho(2c_2\beta(s,0))$. Let $\kappa_0 : \mathbb{R}_{\geq 0} \to \mathbb{R}_{\geq 0}$ be a continuous positive definite function which satisfies $\kappa_0(|\xi|) \leq \rho'(\widetilde{V}(\xi))\frac{1}{2}\kappa(|\xi|)$ for all $\xi \in \mathbb{R}^n$. From (22), we have

$$\begin{aligned}
\nabla V(\xi) \cdot f(\xi,\mu) &\leq -\rho'(\widetilde{V}(\xi))\frac{1}{2}\kappa(|\xi|) + \rho'(\widetilde{V}(\xi))(c_2\gamma_1(|h(\xi)|) + 2c_2\widetilde{\gamma}_2(|\mu|)) \\
&\leq -\kappa_0(|\xi|) + c_2\gamma_1(|h(\xi)|) + 2c_2\widetilde{\gamma}_2(|\mu|) \qquad \forall \xi \in \mathbb{R}^n,\ \forall \mu \in \mathbb{U}.
\end{aligned}$$

This holds at $\xi = 0$ since $V$ is smooth and has a minimum at the origin, so $\nabla V(0) = 0$. Thus $V$ is a smooth iIOSS Lyapunov function for the system. ■